\documentclass[12pt]{amsart}

\newtheorem{pr}{Proposition}[section]
\newtheorem{co}[pr]{Corollary}
\newtheorem{te}[pr]{Theorem}
\newtheorem{lm}[pr]{Lemma}

\theoremstyle{definition}
\newtheorem{de}[pr]{Definition}
\newtheorem{ex}[pr]{Example}
\newtheorem{conj}[pr]{Conjecture}
\newtheorem{que}[pr]{Question}

\theoremstyle{remark}
\newtheorem{re}[pr]{Remark}

\numberwithin{equation}{section}

\def\bea{\begin{eqnarray*}}
\def\eea{\end{eqnarray*}}

\def\toba{{{\mathfrak B}}}
\newcommand{\gr}{\mbox{gr }}
\newcommand{\soc}{\mbox{soc }}

\newcommand{\Z}{\mathbb Z}

\def\phi{\varphi}

\def\qed{{$ \Box $ \vskip 4mm}}
\def\ydh{^{K}_{K}{\mathcal YD}}

\def\Box{\mbox{$\sqcap\!\!\!\!\sqcup$}}

\begin{document}
\title{Co-Frobenius Hopf algebras and the coradical filtration}
\author[Andruskiewitsch and D\u asc\u alescu]{Nicol\'as Andruskiewitsch 
and
Sorin
D\u asc\u alescu}
\address{N. A.: Facultad de Matem\'atica, Astronom\'\i a y F\'\i sica\\
Universidad Nacional de C\'ordoba \\ (5000) Ciudad Universitaria
\\C\'ordoba \\Argentina}
\email{andrus@mate.uncor.edu}
\address{S. D.: Facultatea de Matematica\\
University of Bucharest\\
Str.
Academiei 14\\
RO-70109 Bucharest 1\\
Romania}
\email{sdascal@al.math.unibuc.ro}
\date{January 30, 2001}
\thanks{Parts of this work were done during visits of the second author to the 
University of 
C\'ordoba, Argentina, in May 2000; and of the first author to the University
of Bucharest in November 2000. We thank the FOMEC and 
CNCSIS (Grant C12) for support to these visits.
The first author also thanks ANPCyT, CONICET, CONICOR and Secyt (UNC) for 
partial support. }
\begin{abstract} We prove that a Hopf algebra with a finite coradical 
filtration is co-Frobenius,
{\it i. ~e.} there is a non-zero integral on it. As an application, we show 
that
algebras of functions on quantum groups at roots of one are co-Frobenius. We 
also characterize 
co-Frobenius Hopf algebras with coradical  a Hopf subalgebra. This 
characterization is in the framework of the lifting method, due to H.-J. 
Schneider and the first author.
Here is our main result. Let $H$ be a Hopf algebra whose coradical is
a Hopf algebra. Let $\gr H$ be the associated graded coalgebra and let $R$ be 
the diagram 
of $H$, {\it c.~f.} \cite{as}. Then the following are equivalent: (1) $H$ 
is co-Frobenius; (2) $\gr H$ is co-Frobenius; (3)
$R$  is finite dimensional; (4) the coradical filtration of $H$ is finite.
This Theorem allows to construct systematically examples of co-Frobenius Hopf 
algebras,
and opens the way to the classification of ample classes of such Hopf algebras.
\end{abstract}
\maketitle

\section{Introduction and Preliminaries}

Among the many similarities between the theory of Hopf algebras and 
the theory of groups, the notion of "integral" occupies a central place.
Here, recall that a "left integral" over a Hopf algebra $H$ is a linear map
$\int \in  H^*$ which is left invariant; that is, $\alpha \int = \alpha(1)\int$
for all $\alpha \in H^*$. If $H$ is the algebra of 
regular functions on a compact Lie group $G$, this is exactly (the restriction
of) a left Haar measure 
on $G$. Several basic results on integrals are known, see \cite[Ch. 4]{dnr}.
In particular,  the space of left integrals has dimension $\le 1$. 
A fundamental problem is to determine the class of Hopf 
algebras having a non-zero integral. These are called 
{\it co-Frobenius} Hopf algebras, as explained below.
Classically, finite dimensional Hopf algebras \cite{larsw} and cosemisimple 
Hopf
algebras \cite{sw2} are co-Frobenius. The last class of examples
contains the "compact quantum groups" defined by Woronowicz \cite{woro};
an alternative proof of the existence of a non-zero
Haar measure on a compact quantum group is provided in {\it loc. cit}. 
It is also known that the algebra of 
regular functions on an algebraic group $G$, in characteristic 0,
is co-Frobenius if and only if the group is reductive \cite{sul}.

\medbreak Examples of co-Frobenius Hopf algebras do not abound in the 
literature.
Clearly, the tensor product of two co-Frobenius Hopf algebras also is;
more generally, a cleft extension of two co-Frobenius Hopf algebras also is 
co-Frobenius \cite{bdgn}.
Recently,
several new examples of quantum groups with non-zero integrals have 
been discovered, including  liftings of 
finite quantum linear spaces over an abelian group \cite{bdg}. These are
infinite dimensional 
pointed Hopf algebras which are not cosemisimple. See 
some more examples in \cite{hai}.
One of the corollaries of the main results in this paper is the 
systematic construction of many more new examples of co-Frobenius
Hopf algebras.

\medbreak
The experience of the last 
years shows that an efficient approach to Hopf algebras with 
non-zero integrals is the one from a comodule theory point of view. 
For instance an easy conceptual proof for the uniqueness of the 
integrals was given in this way.

\medbreak

A coalgebra $C$ is called left co-Frobenius if there exists a 
monomorphism of left $C^*$-modules from $C$ to $C^*$. 
This is a generalization of the classical notion of Frobenius algebras.
Also, $C$ is called left semiperfect if the injective envelope $E(S)$ 
of any simple right $C$-comodule $S$ is finite dimensional 
(see \cite{lin}).  
A left co-Frobenius coalgebra is left semiperfect, while the 
converse does not hold in general. However, for a Hopf algebra 
$H$, $H$ is left (or right) co-Frobenius if and only if it 
is left (or right) semiperfect, and this is equivalent to $H$ 
having non-zero integrals. This explains why Hopf algebras 
with non-zero integrals are called co-Frobenius Hopf algebras.

\medbreak

Information about a coalgebra $C$ is captured by its 
coradical filtration $C_0 \subset C_1 \subset \ldots $. 
Since the coradical filtration is a coalgebra filtration, one can construct 
the associated graded coalgebra $\gr C$. Finiteness 
properties of the coalgebras $C$ and $\gr C$ are similar. 
Our first  general result is the following.

\begin{te}\label{A}
Let $C$ be a coalgebra and $\gr C$ the graded coalgebra 
associated to the coradical filtration of $C$. Then $C$ is left 
semiperfect if and only if so is $\gr C$.
\end{te}

\medbreak
Actually, we deduce Theorem  \ref{A} from a more specific result,
relating the injective envelope over $C$ of the
simple right $C$-comodule $S$ with the injective envelope over $\gr C$ of the
simple right $\gr C$-comodule $S$; see Theorem \ref{envgr}.

\medbreak

All the examples we know of co-Frobenius Hopf algebras have a {\it finite}
coradical filtration. The main theme of this paper is:
whether finiteness of the coradical filtration and existence
of non-zero integrals, are equivalent conditions on a Hopf algebra. 
In one direction, a positive answer follows from the following

\begin{te}\label{fincoa}
Let $C$ be a coalgebra such that $C=C_n$ for some $n\geq 0$. Then 
the rational part of the right (or left) $C^*$-module $C^*$ 
is non-zero.
\end{te} 

\begin{co}\label{fincof}
If $H$ is a Hopf algebra such that there exists a non-negative integer $n$ 
with $H_n=H$, then $H$ is co-Frobenius.
\end{co}

We conjecture that the opposite implication is also true:
if $H$ is co-Frobenius, then the coradical filtration is finite.
It is known that the conjecture is true
under the assumption that the coradical is a Hopf subalgebra,
as shown by Radford \cite{rad}.
We present an alternative proof in our main Theorem below. 
Corollary \ref{fincof} allows to present new examples of co-Frobenius
Hopf algebras.
As a distinguished example, we prove  that
algebras of functions on quantum groups at roots of one are co-Frobenius.
The proof uses the "Steinberg-type" decomposition theorem for quantum groups
at roots of one, obtained by Lusztig \cite{lu}.
This material is contained in Section \ref{corfilt}.

\medbreak
An interesting particular situation is when the coradical $H_0$ 
of $H$ is a Hopf subalgebra.Then the coradical filtration of $H$ 
is a Hopf algebra 
filtration, $\gr H$ is a Hopf algebra and there exists 
a Hopf algebra projection $\pi :\gr H\rightarrow H_0$ which splits 
the inclusion of $H_0$ in $\gr H$ as the degree zero component. 
Then the subalgebra $R$ of coinvariants of $\gr H$ with respect 
to the coaction of $H_0$ via $\pi$ has a structure of a Hopf algebra 
in the category $^{H_0}_{H_0}{\mathcal YD}$ of Yetter-Drinfeld modules 
over $H_0$. Moreover $R$, which is called the {\it diagram} of $H$, 
is a graded subalgebra of $\gr H$, and 
$\gr H$ can be reconstructed from $R$ by bosonization, i.e. 
$\gr H$ is isomorphic to the biproduct $R\sharp H_0$. 
This lifting method to study Hopf algebras with the coradical a 
Hopf subalgebra was invented in \cite{as}, and consists in studying first all 
possible 
diagrams $R$, then transferring the information to $\gr H$ by  
bosonization, and finally by lifting the information from 
$\gr H$ to $H$. The method was essentially used only to study 
pointed Hopf algebras, where the coradical is a group algebra. 
In this paper we apply the method to any possible coradical. 
The main result of this paper is 

\medbreak
\begin{te}\label{B} 
Let $H$ be a Hopf algebra with the coradical $H_0$ a Hopf subalgebra. 
Then the following assertions are equivalent.

\begin{enumerate}
\item[(1)]
$H$  is co-Frobenius.

\item[(2)]
 The associated graded Hopf algebra $\gr H$ is co-Frobenius.

\item[(3)]
The diagram $R$ of $H$ is finite dimensional.

\item[(4)]
The coradical filtration $H_0,H_1,\ldots$ of $H$ is finite, i.e. 
there exists $n$ such that $H_n=H$.
\end{enumerate}
\end{te}

\medbreak
This explains for instance why do we need to start with a necessarily finite 
dimensional quantum linear space in order to obtain by the lifting method 
a co-Frobenius Hopf algebra  \cite{bdg}.

The implication (1) $\Rightarrow$ (4)
in Theorem \ref{B}
was proved before in \cite{rad} with different methods.
(1) $\Leftrightarrow$ (2) follows from Theorem \ref{A},
and (4) $\Rightarrow$ (1) is a particular case of Corollary \ref{fincof}.

\medbreak
The coradical of a Hopf algebra is a Hopf subalgebra 
if and only if the tensor product of two simple comodules is completely 
reducible. 
If this holds in a tensor category $\mathcal C$, then it is said
that $\mathcal C$ has the Chevalley property. Furthermore, a
Hopf algebra has the Chevalley property if the category of its finite 
dimensional
modules has the Chevalley property. This denomination honors a classical result
of Chevalley \cite[p. 88]{c}; it was proposed in \cite{aeg}, where triangular
Hopf algebras with the Chevalley property were considered. 
See also \cite{molnar}.

\medbreak

The results we prove in this paper can be used on one hand to construct 
systematically
co-Frobenius Hopf algebras from finite dimensional braided Hopf algebras 
over a cosemisimple Hopf algebra (the coradical), and eventually to classify 
classes
of co-Frobenius Hopf algebras. We perform a preliminary
discussion in Section \ref{examples}.

But on the other hand, these results also serve
as a test to decide that a certain Hopf algebra does not have the 
Chevalley property by looking to the dimension of the injective envelopes 
of simple comodules. See Section \ref{appl}.

\medbreak

We work over a field $k$. The category of right (respectively left) 
comodules over a coalgebra $C$ is denoted by ${\mathcal M}^C$ (respectively 
$^C{\mathcal M}$). If $M$ is a right (or left) $C$-comodule, the injective 
envelope of $M$ in the category ${\mathcal M}^C$ (or  
$^C{\mathcal M}$), which exists since the category of 
comodules is a Grothendieck category, is denoted by $E(M)$. We refer 
to \cite{sw} and \cite{dnr} for notation and facts about coalgebras, 
comodules and Hopf algebras.

\section{The Loewy series of a comodule and the associated graded 
comodule}\label{s1}

Let $C$ be a coalgebra, $C_0\subseteq C_1\subseteq \ldots$  
the coradical filtration of $C$, 
and $M\in {\mathcal M}^C$ a right $C$-comodule with comodule 
structure map $\rho:M\rightarrow M\otimes C$.  
The Loewy series $M_0,M_1,\ldots $ of $M$ is defined as follows: 
$M_0=s(M)$, the socle of $M$, i.e. $M_0$ is the sum of all simple 
subcomodules of $M$; hence for any $n\geq 0$, assume that we  
have defined $M_n$, then $M_{n+1}$ is defined by 
$M_{n+1}/M_n=s(M/M_n)$. We obtain in this way a chain 
$M_0\subseteq M_1\subseteq \ldots $ of $C$-subcomodules of $M$. 
Since $M$ is the sum of all its finite dimensional subcomodules we 
see that $M$ is the union of all $M_n$'s.

\begin{lm}  \label{lema1}
For any $n\geq 0$ we have that $M_n=\rho ^{-1}(M\otimes C_n)$.
\end{lm}

{\bf Proof.} 
Let $J$ be the Jacobson radical of the dual algebra $C^*$, then 
$J=C_0^{\perp}$ 
and $M_n=ann_M(J^{n+1})$ for any $n\geq 0$ (see \cite[Lemma 3.1.9]{dnr}). 
If $m\in M$ such that $\rho (m)\in M\otimes C_n$, then 
since $C_n=(J^{n+1})^{\perp}$, we have that 
$c^*\cdot m= \sum  c^*(m_1)m_0=0$ for any $c^*\in (J^{n+1})$, 
so then $m\in ann_M(J^{n+1})=M_n$. Conversely, let 
$m\in M_n=ann_M(J^{n+1})$. Let $(m_i)_i$ be a basis of $M$ and 
write $\rho (m)=\sum  _im_i\otimes c_i$ for some $c_i\in C$.  
Then for any $c^*\in J^{n+1}$ we have 
$0=c^*\cdot m=\sum  _ic^*(c_i)m_i$, so then $c^*(c_i)=0$ for 
any $i$. Thus $c_i\in (J^{n+1})^{\perp}=C_n$ for any $i$, 
showing that $\rho (m)\in M\otimes C_n$.
\qed

The following result generalizes \cite[Corollary 9.1.7]{sw},  
which states that the coradical filtration of a coalgebra 
is a coalgebra filtration.

\begin{lm} \label{lema2}
For any $n\geq 0$ we have that $\rho (M_n)\subseteq \bigoplus _{i=0}^{n}
M_i\otimes C_{n-i}$.
\end{lm}
{\bf Proof.} We prove by induction on $n$. For $n=0$, we have 
that $M_0$ is the socle of the $C$-comodule $M$, thus it is 
a sum of simple comodules. Since 
the subcoalgebra of $C$ associated to any simple comodule 
is a simple subcoalgebra (see \cite[Exercise 3.1.2]{dnr}), 
we have that $\rho (M_0)\subseteq M_0\otimes C_0$.

\medbreak
Assume now that the assertion holds for $j\leq n$. Let $(u_{0,i})_i$ 
be a basis for $M_0$. We complete this basis with a family 
$(u_{1,i})_i$ up to a basis of $M_1$, then proceed similarly 
by taking families $(u_{2,i})_i,\ldots ,(u_{n+1,i})_i$ for obtaining 
basis of $M_2,\ldots ,M_{n+1}$. Let $m\in M_{n+1}$. Then there 
exist uniquely determined families $(c_{0,i})_i,\ldots ,
(c_{n+1,i})_i$ of elements of $C$ such that 
$$\rho (m)=\sum  _iu_{0,i}\otimes c_{0,i}+\ldots 
+\sum  _iu_{n+1,i}\otimes c_{n+1,i}$$
Denote by $\overline{\rho}:M/M_n\rightarrow M/M_n\otimes C$ the 
comodule structure map of the factor comodule $M/M_n$, and 
by $\overline{x}$ the class modulo $M_n$ of an element $x\in M$. 
Then $\overline{\rho}(\overline{m})=\sum  _i\overline{u_{n+1,i}}
\otimes c_{n+1,i}$. On the other hand $M_{n+1}/M_n$ is a sum 
of simple comodules, so $\overline{\rho}(\overline{m})\in 
M_{n+1}/M_n\otimes C_0$, which shows that $c_{n+1,i}\in C_0$ 
for any $i$. 

\medbreak
On the other hand, since $(\rho \otimes I)\rho =(I\otimes \Delta )\rho$, 
we have that 
$$\sum  _i\rho (u_{0,i})\otimes c_{0,i}+\ldots 
+\sum  _i\rho (u_{n+1,i})\otimes c_{n+1,i}
=\sum  _iu_{0,i}\otimes \Delta (c_{0,i})+\ldots 
+\sum  _iu_{n+1,i}\otimes \Delta (c_{n+1,i})$$
By the induction hypothesis $\rho (u_{r,i})\in 
\bigoplus _{p=0,r}M_p\otimes C_{r-p}$ for any $r\leq n$, so then by looking 
at the terms with the basis element $u_{r,i}$ on the first tensor position,  
where $r\leq n$, we see that 
$$\Delta (c_{r,i})\in C_0\otimes C+C_1\otimes C+\ldots +C_{n-r}\otimes C
+C\otimes C_0=C_{n-r}\otimes C+C\otimes C_0$$ 
This implies that $c_{r,i}\in C_{n-r}\wedge C_0=C_{n-r+1}$, which ends 
the proof.
\qed

Now we construct the graded comodule associated to the Loewy series  
of a $C$-comodule $M$. In particular, for $M=C$, we obtain the 
graded coalgebra associated to the coradical filtration of $C$. 
For any $i\geq 1$ we denote by $\pi _i:C\rightarrow C/C_{i-1}$ 
and $p_i:M\rightarrow M/M_{i-1}$ the natural projections. 
We consider the spaces $\gr C=\oplus _{i\geq 0}C_i/C_{i-1}$ 
and $\gr M=\oplus _{i\geq 0}M_i/M_{i-1}$, where we take  
$C_{-1}=0, M_{-1}=0$. By Lemma \ref{lema2}, the map 
$$(\oplus _{i=0}^{n}p_i\otimes \pi _{n-i})\circ \rho :M\rightarrow 
\oplus _{i=0}^{n}(M/M_{i-1}\otimes C/C_{n-i-1})$$
induces a linear map 
$${\rho}_n :M_n\rightarrow 
\oplus _{i=0}^{n}(M_i/M_{i-1}\otimes C_{n-i}/C_{n-i-1})$$
Using again Lemma \ref{lema2} we have that 
$\rho _n(M_{i-1})=0$, and thus $\rho _n$ induces a linear 
map 
$$\overline{\rho}_n :M_n/M_{n-1}\rightarrow 
\oplus _{i=0}^{n}(M_i/M_{i-1}\otimes C_{n-i}/C_{n-i-1})$$
If we regard $\oplus _{i=0}^{n}(M_i/M_{i-1}\otimes C_{n-i}/C_{n-i-1})$ 
as a subspace of $\gr M\otimes \gr C$, the sum of all 
$\overline{\rho}_n$'s define a linear map 
$\overline{\rho}:\gr M\rightarrow \gr M\otimes \gr C$. 
This map can be described as follows. For $m\in M_n$, write 
$\rho (m)=\sum  _{i=0}^{n}m_{0,i}\otimes m_{1,n-i}$, 
a Sigma Notation type representation, with the $m_{0,i}$'s lying 
in $M_i$ and the $m_{1,n-i}$'s lying in $C_{n-i}$.
Then we have that $\overline{\rho}(p_n(m))=
\sum  p_i(m_{0,i})\otimes \pi _{n-i}(m_{1,n-i})$; and this does not depend on 
the chosen representation of $\rho (m)$. In the case where $M=C$ and 
$\rho=\Delta$, we obtain the map 
$\overline{\Delta}:\gr C\rightarrow \gr C\otimes \gr C$, 
defined by $\overline{\Delta}(\pi _n(c))=
\sum  \pi _i(c_{1,i})\otimes \pi _{n-i}(c_{1,n-i})$ for $c\in C_n$ and  
any representation $\Delta (c)=\sum  c_{1,i}\otimes c_{2,n-i}$ with the 
same convention as above. Straightforward computations show that 
$\gr C$ is a graded coalgebra with comultiplication $\overline{\Delta}$ 
and $\gr M$ is a graded right $\gr C$-comodule via $\overline{\rho}$.
We denote by $\gr M(n)=M_n/M_{n-1}$ the homogeneous component of 
degree $n$ of $\gr M$. 

The construction of the graded coalgebra associated to a coalgebra filtration 
of a coalgebra goes back to Sweedler's book (see \cite[Section 11.1]{sw}).  
By \cite[Lemma 2.3]{as} we have that $\gr C$ is coradically graded, 
i.e. its coradical filtration is given by $(\gr C)_n=
\oplus _{i\leq n}\gr C(i)$ for any $n\geq 0$.

\medbreak
\begin{te}\label{envgr} Let $S$ be a simple right $C$-comodule, $C$ any 
coalgebra.
The graded comodule associated to the injective envelope over $C$ of
 $S$, is isomorphic to the the injective envelope over $\gr C$ of
the simple right $\gr C$-comodule $S$.
\end{te}
{\bf Proof.}
Since the socle of a direct sum of comodules is the sum of the socles of these 
comodules, we see that for any family $(M(\lambda))_{\lambda \in \Lambda}$ 
of right $C$-comodules we have $\gr (\oplus _{\lambda \in \Lambda}
M(\lambda ))$ $\simeq \oplus _{\lambda \in \Lambda}\gr M(\lambda )$ 
as right $\gr C$-comodules. In particular 
if we write $C_0=Soc(C)=\oplus _{S\in {\mathcal S}} \, S$ for some family 
$\mathcal S$ of 
simple right $C$-subcomodules of $C$, then $C=\oplus _{S\in {\mathcal S}} \, 
E(S)$ 
(see \cite[Theorem 2.4.16]{dnr}), so then 
$\gr C=\oplus _{S\in {\mathcal S}} \, \gr E(S)$ 
as right $\gr C$-comodules. In particular $\gr E(S)$ is injective 
as a right $\gr C$-comodule. Since $S=E(S)_0=
\gr E(S)(0)\subseteq \gr E(S)$, we obtain that 
$E_{\gr} (S)\subseteq \gr E(S)$, where 
$E_{\gr} (S)$ is an injective envelope of $S$ 
as a right $\gr C$-comodule. 

On the other hand $(\gr C)_0=\gr C(0)=C_0=
\oplus _{S\in {\mathcal S}} \, S$ as right $\gr C$-comodules, so then  
$\gr C=\oplus _{S\in {\mathcal S}} \, E_{\gr} (S)$, showing that 
$E_{\gr} (S)=\gr E(S)$ for any $S\in {\mathcal S}$. 
\qed

{\bf Proof of Theorem  \ref{A}.} Given a simple $S$, $E(S)$  is finite 
dimensional if and only if 
so is $\gr E(S)\simeq E_{\gr} (S)$.
\qed

\section{Coradical filtration of a co-Frobenius Hopf algebra}\label{corfilt}

Let $H$ be a Hopf algebra and $H_0, H_1,\ldots$ its coradical filtration. 
We first show that if the coradical filtration is finite, then 
necessarily $H$ is co-Frobenius. This follows from  
Theorem \ref{fincoa}, which gives information about coalgebras with 
finite coradical filtration.

{\bf Proof of Theorem \ref{fincoa}.} If $C=C_0$, then $C$ is cosemisimple, 
therefore it is left and right  
co-Frobenius (see \cite[Exercise 3.3.17]{dnr}). In particular 
the rational part of $C^*$ as a right (or left) $C^*$-module 
is non-zero. 
If $C$ is not cosemisimple, let $n$ be such that 
$C_{n-1}\neq C_n=C$. Then $C/C_{n-1}$ is a semisimple right $C$-comodule, 
in particular it contains a maximal subcomodule. This induces a 
maximal right $C$-subcomodule $X$ of $C$. 
The proof goes now as in \cite[Proposition 3.2.2]{dnr}. 
The natural projection 
$C\longrightarrow C/X$ produces an injective morphism of right  
$C^*$-modules $(C/X)^*\hookrightarrow C^*$. But $C/X$ is simple, 
so it is finite dimensional and rational as a left $C^*$-module, 
and then $(C/X)^*$ is a rational right $C^*$-module (see 
\cite[Lemma 2.2.12]{dnr}). This shows that the rational part of the 
right $C^*$-module $C^*$ is non-zero.  
\qed

Corollary \ref{fincof} follows from Theorem \ref{fincoa} and the 
fact that a Hopf algebra $H$ is co-Frobenius if and only if 
$H^{*rat}\neq 0$ \cite[Ch. 5]{dnr}. 
Note that for a Hopf algebra the left and the right 
rational parts of $H^*$ are equal, and they are denoted by 
$H^{*rat}$.

As said in the Introduction, we conjecture that the converse 
of Corollary \ref{fincof} is true: 

\begin{conj}\label{converse}The coradical filtration of a  co-Frobenius Hopf 
algebra is finite. 
\end{conj}

We would like to mention the 
following more precise question, which arose in joint work of 
Sonia Natale and the first-named author.

\begin{que}\label{que1} Let $H$ be a Hopf algebra and let $T$ be a left 
integral on $H$. 
If $H_m$ is different from $H$, is it true that $T$ vanishes 
on $H_m$?\end{que}

Note that a positive answer to Question \ref{que1} implies  
Conjecture \ref{converse}. For $m=0$, the answer 
is positive by Maschke's Theorem for Hopf algebras \cite{sw2}.

We now show how Corollary \ref{fincof} allows to determine that some 
important Hopf algebras are co-Frobenius. We begin by the following 
general Lemma.

\begin{lm}\label{lemgral}
Let $H$ be a  Hopf algebra and let $K$ be a Hopf
subalgebra such that
\begin{flalign}
\label{5.1}
& KH_{0} \subseteq H_{0} \text{ (in particular, $K \subseteq H_{0}$).}& 
\\ \label{5.2}
& H  \text{ is of finite type as left module via multiplication over } K. &
\end{flalign}
Then the coradical filtration of $H$ is finite; and $H$ is co-Frobenius.
\end{lm}

{\bf Proof.} By \eqref{5.1}, each term $H_n$ of the coradical filtration is a 
$K$-submodule; 
by \eqref{5.2}, the coradical filtration is then finite. The last claim 
follows from
Corollary \ref{fincof}. \qed

In representation-theoretic terms, condition \eqref{5.1} means the 
following. If
$V$ is a simple $K$-comodule and $W$ is a simple $H$-comodule then
the $H$-comodule $V\otimes W$ is completely reducible. 
Indeed, let  $\widehat{H}$ denote the set of (isomorphism classes of)
simple comodules over a Hopf algebra $H$; we will confuss a class with a 
representant
without danger. If $U$ is any finite dimensional $H$-comodule, 
then let $C_U$ denote the space of matrix coefficients of $U$;
it is a subcoalgebra of $H$. We have $C_{U \otimes W} = C_UC_W$.
Also $K = \sum_{V\in \widehat{K}} C_V$, $H_0 = \sum_{W\in \widehat{H}} C_W$, 
thus the claim.

We are led to the following definition.

\begin{de} Let $H$ be a Hopf algebra with bijective antipode.
Let $\widehat{H}_{\soc}$ be the subset of $\widehat{H}$ consisting of all
simple comodules $V$ such that
\begin{equation} V\otimes W \text{ and } W\otimes V \text{ are completely 
reducible }
\end{equation}
for all $W\in \widehat{H}$. The Hopf socle of $H$ is
$$
{H}_{\soc} = \sum_{V\in \widehat{H}_{\soc}} C_V.
$$
\end{de}

\begin{lm}\label{socle} Let $H$ be a Hopf algebra with bijective antipode.
The Hopf socle of $H$ is a cosemisimple Hopf subalgebra of $H$.\end{lm}
{\bf Proof.} It is clear that ${H}_{\soc}$ is a subcoalgebra of $H_0$.
We  prove that ${H}_{\soc}$ is a subalgebra of $H$.
Let $V, U \in \widehat{H}_{\soc}$. Then $V\otimes U = \oplus_{1\le j\le r} 
W_j$,
where $W_j$ are simple comodules.
We have to show that $W_j \in \widehat{H}_{\soc}$, $1\le j\le r$. Let $W\in 
\widehat{H}$.
Then $V\otimes U \otimes W$ is completely reducible; since $W_j \otimes W$
is a subcomodule of $V\otimes U \otimes W$, it is also completely reducible.
By an analogous argument, $W \otimes W_j$ is completely reducible.
Therefore, $W_j \in \widehat{H}_{\soc}$. We finally prove that ${H}_{\soc}$
is  stable under the antipode $\mathcal S$. If $V$ is a finite dimensional 
comodule,
we denote by $V^*$ be the comodule structure on the dual  defined via 
$\mathcal S$,
and by $^*V$, the comodule structure on the dual  defined via $\mathcal 
S^{-1}$.
We know that $(V\otimes W)^* \simeq W^* \otimes V^*$, $^*(V\otimes W) \simeq 
{}^*W \otimes {}^*V$
for any two finite dimensional comodules $V$ and $W$. If $V \in 
\widehat{H}_{\soc}$
and $W\in \widehat{H}$, then
$^*(V^*\otimes W) \simeq {}^*W \otimes V$ is completely reducible, therefore
$V^*\otimes W$ is completely reducible. Similarly, 
$W\otimes V^*$ is completely reducible. Thus, $V^* \in \widehat{H}_{\soc}$,
as needed. \qed

Lemma \ref{lemgral} immediately implies:

\begin{co}\label{cor-socle}Let $H$ be a Hopf algebra with bijective antipode.
If $H$ is of finite type as left module via multiplication over ${H}_{\soc}$,
then $H$ is co-Frobenius. \qed\end{co}

We ignore if the converse of Corollary \ref{cor-socle} is true;
that is, if any co-Frobenius Hopf algebra is of finite type over ${H}_{\soc}$.

\medbreak
From now on, we assume that $k$ is algebraically closed and has characteristic 
0.
Let $\mathfrak g$ be a simple finite dimensional Lie algebra of rank $n$, 
and let $G$ be the
corresponding simply-connected algebraic group. Let $q\in k$ be a root of 1 of 
odd order 
$\ell \ge 3$, $\ell$ not divisible by 3 if  $\mathfrak g$ contains a component 
of type $G_2$; and let $U_q(\mathfrak g)$ be the quantized enveloping algebra 
as defined in 
\cite{lu}. Recall the definition of modules of type 1 \cite{lu}. 
The following facts were proved by Lusztig:

\begin{itemize}
\item The set of isomorphism classes of simple $U_q(\mathfrak g)$-modules of 
type 1 is parametrized
by $(\Z_{\ge 0})^n$ \cite[Prop. 6.4]{lu}. If $\lambda \in (\Z_{\ge 0})^n$,
the corresponding simple module $L(\lambda)$ has highest weight $\lambda$.

\medbreak

\medbreak
\item If  $\lambda \in (\Z_{\ge 0})^n$, decompose it as  $\lambda = \lambda' +
\ell \lambda''$, where $\lambda', \lambda'' \in (\Z_{\ge 0})^n$ with each 
entry of
$\lambda'$ living in the interval $[0, \ell -1]$. Then $L(\lambda) \simeq 
L(\lambda')
\otimes L(\ell\lambda'')$ (Steinberg-type Theorem, \cite[Th. 7.4]{lu}).

\medbreak
\item There is an epimorphism of Hopf algebras $U_q(\mathfrak g) \to 
U(\mathfrak g)$. 
If $\lambda'' \in (\Z_{\ge 0})^n$, the simple $U_q(\mathfrak g)$-module  
$L(\ell\lambda'')$
inherits a structure of $U(\mathfrak g)$; and as such, it is simple with 
highest weight $\lambda''$
\cite[Prop. 7.5]{lu}. Denote this last module by ${\underline L}(\lambda'')$.
\end{itemize}

Let $k_q[G] \subset U_q(\mathfrak g)^*$ be the Hopf algebra of matrix 
coefficients
of modules of type 1. This is the algebra of functions on a quantum group at a 
root of one, as
in \cite{dcl}.  The coradical of $k_q[G]$ is the span of the matrix 
coefficients of 
the simple modules $L(\lambda)$, $\lambda \in (\Z_{\ge 0})^n$. 
It is known that it is {\it not} a Hopf subalgebra.

The subalgebra of $k_q[G]$ spannned by the matrix coefficients of 
the simple modules $L(\ell\lambda'')$, $\lambda'' \in (\Z_{\ge 0})^n$, is
a Hopf subalgebra of $k_q[G]$, isomorphic to the
algebra $k[G]$ of regular functions on $G$. 
It is known that it is central \cite{dcl}.

\begin{te}\label{quantumfunctalgs}
$k_q[G]$ is co-Frobenius. \end{te}

{\bf Proof.} Let $H = k_q[G]$, $K = k[G]$. \eqref{5.2} holds by \cite{dcl}.
Condition \eqref{5.1} is a consequence of the "Steinberg-type" decomposition 
theorem evoked above. Indeed, let $V = L(\ell\mu)$, $W = L(\lambda)$
with  $\mu, \lambda \in (\Z_{\ge 0})^n$. Decompose $\lambda$ as before,
so that
$L(\lambda)\simeq L(\lambda')
\otimes L(\ell\lambda'')$, $\lambda', \lambda'' \in (\Z_{\ge 0})^n$. 
By Weyl's Theorem,  ${\underline L}(\lambda'') \otimes {\underline L}(\mu)$
is completely reducible, say isomorphic to $\oplus_{j=1, \dots, M} {\underline 
L}(\sigma_j)$.
Then
\begin{align*}
W\otimes V & \simeq L(\lambda') \otimes L(\ell\lambda'') \otimes L(\ell\mu)
\simeq L(\lambda') \otimes \left( \oplus_{j=1, \dots, M} 
L(\ell\sigma_j)\right) \\
& \simeq \oplus_{j=1, \dots, M} L(\lambda') \otimes L(\ell\sigma_j) 
\simeq \oplus_{j=1, \dots, M} L(\lambda' + \ell\sigma_j).
\end{align*}
Therefore, $KH_{0} = H_{0}K \subseteq H_{0}$. 
The Theorem  follows now from Lemma \ref{lemgral}. \qed

\begin{re} It is known that $k_q[G]$ is an extension of $k[G]$ by a finite 
Hopf algebra,
namely the dual of the corresponding Frobenius-Lusztig kernel. However,
we do not known if it is a {\it cleft} extension; otherwise we could
have deduced Theorem \ref{quantumfunctalgs} from \cite[5.2]{bdgn}.
It is however known that $k_q[G]$ is free over $k[G]$ \cite{bg, bgs}. \end{re}

\section{Co-Frobenius Hopf algebras whose coradical is 
a Hopf subalgebra}\label{cofrobcorad}

In this section we assume that the coradical $H_0$ of $H$ 
is a Hopf subalgebra.  

Write as in Section \ref{s1} $H_0=
\oplus _{S\in {\mathcal S}} \, S$ for some family $\mathcal S$ of 
simple left $H$-subcomodules of $H$. We identify   
$\gr H$ and $R\sharp H_0$, then 
$\gr H=\oplus _{S\in {\mathcal S}} \, R\sharp S$. 

\begin{lm}\label{2.2}
$R\sharp S$ is an injective left $\gr H$-subcomodule of $\gr H$ for 
any $S\in {\mathcal S}$.
\end{lm}

{\bf Proof.} The comultiplication of the biproduct is 
$\Delta (r\sharp h)=\sum  r^1\sharp (r^2)_{-1}h_1\otimes (r^2)_0\sharp h_2$,  
where $r\mapsto \sum  r^1\otimes r^2$ is the comultiplication of the braided 
Hopf 
algebra $R$, and $r\mapsto \sum  r_{-1}\otimes r_0$ is the left coaction of 
$H_0$ 
on $R$. This implies that $\Delta (R\sharp S)\subseteq \gr H\otimes 
(R\sharp S)$. The injectivity follows from the fact that $R\sharp S$ is a 
direct 
summand of the left $\gr H$-comodule $\gr H$. 
\qed

\begin{co}    \label{anvinj}
Let $S\in {\mathcal S}$ and let $E_{\gr} (S)$ be the injective envelope of 
$S$ as a left $\gr H$-comodule. Then $E_{\gr} (S)=R\sharp S$.
\end{co}

{\bf Proof.} 
The result follows immediately from the previous lemma if 
we take into account that $S$ is contained in $R\sharp S$, 
and that for any $S,T\in {\mathcal S}$, $S\neq T$,  $(R\sharp S)\cap (R\sharp 
T)=0$.
\qed

{\bf Proof of Theorem \ref{B}} $(1)\Leftrightarrow (2)$ follows from 
Theorem  \ref{A}. 

\medbreak
$(2)\Leftrightarrow (3)$ follows from Corollary \ref{anvinj}.

\medbreak  
$(3)\Rightarrow (4)$. Since $R$ is finite dimensional there exists a positive  
  
integer $n$ such that $R=R(0)+\ldots +R(n)$. Then 
$\gr H=\gr H(0)+\ldots +\gr H(n)$, implying that $H_n=H$.

\medbreak
$(4)\Rightarrow (1)$ follows from Corollary \ref{fincof}.
\qed

Here is an alternative proof of $(4)\Rightarrow (1)$. If the coradical 
filtration is finite, then the graded braided Hopf algebra $R$ has a finite 
grading, {\it i.e.}
$R = \oplus_{0\le n \le N } R(n)$, with $R(N) \neq 0$. Since $R(0) = k$, any 
element 
$t\in R(N)$, $t\neq 0$ is an integral in $R$, see for instance \cite[Prop. 
3.2.2]{AG}.
Hence $R$ is finite dimensional by \cite{FMS}. 

We state another consequence of Corollary \ref{anvinj}. 
Let $E$ be the injective envelope of the trivial  $H$-comodule $k1$.

\begin{co}    \label{otrocoro}
Let $S\in {\mathcal S}$ and let $E (S)$ be the injective envelope of 
$S$ as a left $H$-comodule. Then $E(S) \simeq E \otimes S$.
\end{co}

{\bf Proof.} Observe first that the claim is true for $\gr H$, by Corollary 
\ref{anvinj}.
 Since $E \otimes S$ is a direct summand of $H \otimes S \simeq H^{\dim S}$,
it is injective. Thus there exists a monomorphism of $H$-comodules 
$\phi: E(S) \to E \otimes S$. As $H_0$ is a Hopf subalgebra of $H$, 
the Loewy series of $E \otimes S$ is $(E \otimes S)_n = E_n \otimes S$. 
Therefore,
$\gr \phi: \gr E(S) \to \gr (E \otimes S) \simeq (\gr E) \otimes S$ is an 
isomorphism.
By a standard argument, $\phi$ is an isomorphism. \qed

\section{Some applications}\label{appl}

We give now some numerical criteria for the Chevalley property.
All of them follow from Lemma \ref{2.2}.

\begin{pr}
Let $H$ be a finite dimensional Hopf algebra with 
a simple right (or left) $H$-comodule $S$ such that 
the dimension of $S$ does not divide the dimension of $E(S)$. 
Then the category of finite dimensional $H$-comodules 
does not have the Chevalley property. \qed
\end{pr}
If the coradical of a Hopf algebra $H$ is a Hopf subalgebra, then
the diagram $R$ is isomorphic to the injective envelope of the trivial $\gr  
H$-comodule. 
Then:

\begin{pr}
Let $H$ be a finite dimensional Hopf algebra such that 
the dimension of the injective envelope of the trivial $H$-comodule does not 
divide the dimension of $H$. 
Then the category of finite dimensional $H$-comodules 
does not have the Chevalley property. \qed
\end{pr}

Dualizing, we have:

\begin{pr}
Let $H$ be a finite dimensional Hopf algebra with 
a simple right (or left) $H$-module $S$ such that 
the dimension of $S$ does not divide the dimension of the projective cover 
$P(S)$. 
Then  $H$ does not have the Chevalley property. \qed
\end{pr}

\begin{pr}
Let $H$ be a finite dimensional Hopf algebra such that 
the dimension of the projective cover of the trivial $H$-module 
does not divide the dimension of $H$.  
Then $H$ does not have the Chevalley property. \qed
\end{pr}

If $C$ is a coalgebra and $M\in {\mathcal M}^C$ a right $C$-comodule, 
we can consider the Poincar\'e series  of $M$, namely
$$\ell(M) = \sum_{n\ge 0} \ell_n(M) X^n \in {\mathbb Z}[[X]], \qquad 
\text{where} 
\quad \ell_n(M) = \dim M_{n}/M_{n - 1}. $$ 
Clearly, $\ell(M) = \ell(\gr M)$.

\begin{pr}
Let $H$ be a co-Frobenius  Hopf algebra with the coradical a Hopf subalgebra.
If $S$ is a simple right comodule, then the Poincar\'e polynomial of the 
injective envelope
$E(S)$ satisfies a "Poincar\'e duality": $\ell_n(E(S)) = \ell_{\text{top} - 
n}(E(S))$,
for all $n$, $0\le n \le \text{top}$. Here "top" means the degree of the 
Poincar\'e polynomial.
\end{pr}

{\bf Proof.} By Corollary \ref{anvinj}, we have
$$
\ell(E(S)) = \ell(\gr E(S)) = \ell(R \# S) = \ell(R) \dim S.
$$
Then the claim follows from the analogous claim for $R$, which is well-known 
\cite{ni}.
\qed

\section{Examples from the Lifting Method}\label{examples}

We now explain how to obtain new examples, and eventually
classification results, of co-Frobenius Hopf algebras via Theorem \ref{B}.
We rely on the Lifting method \cite{as, as2}; a detailed exposition is
\cite{as3}. We fix a cosemisimple Hopf algebra $K$ and seek for 
co-Frobenius Hopf algebras $H$ with $H_0 \simeq K$.

\medbreak
Recall that a braided vector space is a pair $(V, c)$,
where $V$ is a vector space (which we assume finite dimensional)
and $c: V\otimes V \to V\otimes V$ is an isomorphism satisfying the braid 
equation.
Given a braided vector space $(V, c)$, there is a remarkable 
graded braided Hopf algebra ${\toba} (V)$ called the Nichols algebra of $(V, 
c)$. 
The Lifting method in the present setting consists of the following steps.

\medbreak 
(a).  Determine when ${\toba}(V)$  is finite dimensional, for all braided 
vector spaces $(V,c)$ arising as Yetter-Drinfeld modules over {\it some} 
cosemisimple Hopf algebra.

\medbreak 
(b).  For those $V$ as in (a), find in how many ways, if any, they can be 
actually realized as Yetter-Drinfeld modules over our fixed $K$.

\medbreak 
(c). For ${\toba}(V)$ as in (a), compute all Hopf algebras $H$ such that $\gr 
H \simeq 
{\toba}(V) \# K$ ("lifting").

\medbreak 
(d). Investigate whether any {\it finite dimensional} graded braided Hopf 
algebra $R = \oplus_{n\ge 0} R(n)$ in $\ydh$ 
satisfying $R(0) = k 1$ and $P(R) = R(1)$,  is generated by its primitive 
elements, {\it i. e.} is a Nichols algebra.

\medbreak For the discussions of the different steps, 
we can profit previous investigations
in the setting of finite dimensional Hopf algebras \cite{as, as2}.

\medbreak Let us first discuss step (a).
Actually, we have an ample supply of finite dimensional Nichols
algebras ${\toba} (V)$.  Let us recall some definitions from \cite{as, as2}.

\begin{itemize}
\item A braided vector space
$(V, c)$ is of {\it diagonal type} if $V$ has a basis
$x_{1}, \dots, x_{\theta}$ such that $c(x_{i} \otimes x_{j})
= q_{ij} x_{j} \otimes x_{i}$, for all $1\le i, j\le \theta$,
where $q_{ij}$ are some scalars.
\end{itemize}
Since we are 
interested in finite dimensional Nichols algebras,
we can (and shall) assume that $q_{ii} \neq 1$, for all $i$.

\begin{itemize}
\item A braided vector space of diagonal type is
of {\it Cartan type} if $q_{ij}q_{ji}
= q_{ii}^{a_{ij}}$, for all $1\le i \neq j\le \theta$,
where $a_{ij}$ are some non-positive integers \cite{as2}. Set
$a_{ii} = 2$ for all $i$.
\end{itemize}

The integers $a_{ij}$ can be chosen so that the matrix
$(a_{ij})$ is a generalized Cartan matrix.  The following Theorem was proved
in \cite{as2}, from results of Lusztig and using the twisting operation.

\begin{itemize}
\item If $(a_{ij})$ is actually a finite Cartan matrix
and the orders of the $q_{ij}$'s are odd and not divisible by 3, when
the matrix  $(a_{ij})$ has a component of type $G_2$, then ${\toba}(V)$
has finite dimension (which can be explicitly computed). \end{itemize}

\medbreak The simplest example is when $q_{ij}q_{ji}
= 1$ for all $i\neq j$. Then $\toba(V)$ is called a finite quantum linear 
space;
the examples in \cite{bdg} are exactly liftings of quantum linear spaces over 
an arbitrary abelian
group.

\medbreak 
Let us comment on step (b) in the setting of the examples just discussed.
Let us fix a braided vector space $(V, c)$
of finite Cartan type, with the restrictions  
above. To realize $(V, c)$ in $\ydh$, for our fixed cosemisimple Hopf algebra 
$K$,
we need:

\begin{itemize}
\item a family of characters $\chi_1, \dots, \chi_{\theta}$ of $K$,
\item and a family of central group-likes $g_1, \dots, g_{\theta} \in G(K) 
\cap Z(K)$
such that  \end{itemize}
$$\chi_j(g_i) = q_{ij}, \quad 1\le i, j\le \theta.$$

For instance, $(V, c)$ can be realized as a Yetter-Drinfeld module over 
the group algebra of ${\mathbb Z}^s$, if $s\ge \theta$. 
To determine whether it can be realized as a Yetter-Drinfeld module over 
the group algebra of ${\mathbb Z}^s$, for $s < \theta$ is an interesting 
arithmetical problem.

\begin{ex} Let $s = 1$. 
If all the orders of the entries $q_{ij}$ divide a fixed odd prime number $p > 
3$,
then necessarily $\theta \le 2$ and the following restrictions on $p$ are in 
force:

\begin{itemize}
\item If the Cartan matrix is of type $A_2$, then $p \equiv 1 \mod 3$.

\smallbreak
\item If the Cartan matrix is of type $B_2$, then $p \equiv 1 \mod 4$.

\smallbreak
\item If the Cartan matrix is of type $G_2$, then $p \equiv 1 \mod 
3$.\end{itemize}

Indeed, the image of ${\mathbb Z}$ in Aut $V$ is a cyclic group of order $p$,
and we can apply \cite[Th. 1.3]{as2}. A similar discussion is also valid for 
$p = 3$,
see {\it loc. cit}. We obtain in this way many new examples of co-Frobenius 
Hopf algebras,
namely ${\toba}(V) \# k{\mathbb Z}$.
The situation is somewhat different for $p=2$ \cite{ni, cd}; the examples 
arising here
are particular cases of those in \cite{bdg}.
\end{ex}

\medbreak
The computation of the liftings of the Hopf algebras ${\toba} (V) \# K$,
step (c) of the method,
requires the knowledge of a presentation by generators and relations of 
${\toba} (V)$. This was obtained in  \cite{as4} for braidings of finite Cartan 
type.
It is likely that the techniques of the finite dimensional case are also
useful here, see \cite{as4, as3}. Finally, let us mention concerning step (d),
that a positive result in this direction is given in \cite{as4}.

\medbreak
Let us come back to the case $K = k{\mathbb Z}$. It is natural to ask whether
there are Yetter-Drinfeld modules $V$ over $K$ which are {\it not} of finite 
Cartan type
but such that ${\toba} (V)$ is finite dimensional. First, it can be shown 
that such $V$ is necessarily of diagonal type. Second, a small quantity
of examples of such $V$ are known \cite{ni, gr}.

\medbreak Are there examples of Yetter-Drinfeld modules
$V$ over some cosemisimple Hopf algebra $K$ which are {\it not} of diagonal 
type,
but such that ${\toba} (V)$ is finite dimensional? Yes, a few ones; see 
\cite{MiS, gr}. They can be realized over suitable group algebras.

\medbreak On another direction, let $K$ be the algebra of regular functions on 
a
simple algebraic group. Then we do not know any example of Yetter-Drinfeld 
module 
$V$ over $K$ with ${\toba} (V)$ finite dimensional.

\bigbreak

\end{document}